\documentclass[11pt]{amsart}
\usepackage{amsmath,amsthm,epsf,amscd,amssymb,verbatim}

\newtheorem{thm}{Theorem}[section]
\newtheorem{lem}[thm]{Lemma}

\newtheorem{prop}[thm]{Proposition}
\newtheorem{conj}[thm]{Conjecture}

\theoremstyle{definition}
\newtheorem{dfn}[thm]{Definition}

\theoremstyle{remark}
\newtheorem{remark}[thm]{Remark}

\numberwithin{equation}{section}

\newcommand{\bbr}{\begin{remark}}        
\newcommand{\eer}{\end{remark}}

\newcommand{\df}[1]{{\em{#1}}}

\newcommand{\tm}[1]{Theorem~{\ref{#1}}}
\newcommand{\eqn}[1]{{Equation~\ref{#1}}}

\font\bbb=msbm10 scaled 1100

\newcommand{\be}{\begin{equation}}
\newcommand{\ee}{\end{equation}}
\newcommand{\bea}{\begin{eqnarray}}
\newcommand{\eea}{\end{eqnarray}}
\newcommand{\bmini}{\footnotesize\begin{center}\begin{minipage}{5.5in}}
\newcommand{\emini}{\end{minipage}\end{center}\normalsize}
\newcommand{\pf}{{\em Proof: }}

\newcommand{\R}{\mbox{\bbb R}}
\newcommand{\C}{\mbox{\bbb C}}
\newcommand{\real}{\mbox{\bbb R}}       
\newcommand{\eg}{{\em e.g.}}
\newcommand{\ie}{{\em i.e.}}

\newcommand{\mathspace}{\;\;\;\;\;}
\newcommand{\del}{\partial}

\begin{document}

\title{Contact topology and hydrodynamics II: solid tori}
\author{John Etnyre}
\address{Department of Mathematics, Stanford University, 
	Palo Alto, CA, 94305}
\thanks{JE supported in part by NSF Grant \# DMS-9705949.}

\author{Robert Ghrist}
\address{School of Mathematics, Georgia Institute of Technology,
	Atlanta, GA 30332}
\thanks{RG supported in part by NSF Grant \# DMS-9971629.}


\begin{abstract}
We prove the existence of periodic orbits for steady $C^\omega$
Euler flows on all Riemannian solid tori. By using the correspondence
theorem from part I of this series, we reduce the problem to the Weinstein 
Conjecture for solid tori. We prove the Weinstein Conjecture on the
solid torus via a combination of results due to Hofer et al. and 
a careful analysis of tight contact structures on solid tori.
\end{abstract}

\maketitle

\section{Introduction and summary}

Problems associated with the existence of periodic orbits in 
flows differ sharply from their discrete counterparts. The 
index-theoretic methods available for detecting periodic orbits
of two-dimensional self-diffeomorphisms are not sufficient for 
understanding general three-dimensional flows. 
Indeed, the recent progress on the Seifert Conjecture by 
K. Kuperberg and G. Kuperberg \cite{Kup94,KK94} has made it clear
that, for flows of arbitrary 
regularity, there is ``too much room'' to have a topological 
forcing theory: additional constraints are required, many of which
should be geometric in nature.
The interesting problem is now to find sharp boundaries 
on the space of vector fields in dimension three which 
separate those fields without periodic orbits. Currently, there
is great interest in the case of volume-preserving and Hamiltonain 
flows, since, by a classical theorem of Poincar\'e, 
almost all orbits are recurrent.
The Kuperberg plug constructions do not work in this 
category (but see \cite{Kup96b} for a $C^1$ construction). 

In this series of papers, we are concerned with periodic orbits in 
a particular class of volume-preserving vector fields which models 
the motion of the simplest possible fluid and plasma flows: these are 
the steady, perfect, incompressible fields, or {\em Euler fields}. 
In the realm of fluid and plasma dynamics, periodic orbits 
play a naturally important role. For example, the existence of 
a hyperbolic periodic orbit in a steady Euler flow is sufficient
to conclude hydrodynamic instability of the solution \cite{FV92}.
There are several connections 
between the embedding properties of periodic orbits and physical 
properties of fluids/plasmas, such as energy bounds, helicity, 
and the possibility of finite-time singularities (see, \eg, 
\cite{Mof85,Mof94}). In part III of this series \cite{EG:III} we 
consider the knot theory of periodic orbits in Euler flows.

In part I of this series \cite{EG:I}, we proved the existence of periodic 
solutions to all $C^\omega$ steady nonsingular solutions to the Euler 
equations for an inviscid fluid flow on a Riemannian 3-sphere. In 
applications, fluid flows on $S^3$ are not terribly significant. 
In \cite{EG:I}, algebraic-topological conditions were also derived 
under which steady nonsingular Euler flows on a Riemannian cube 
$[0,1]^3$ with periodic boundary conditions possess a closed orbit.
While this global geometry is a common domain in theoretical 
fluid dynamics (Fourier analysis is particularly simple here), 
one cannot reasonably call $T^3$ ``physical.'' Perhaps the most 
interesting case for applications (\eg, tokamaks) involves compact 
domains in $\real^3$. Under the condition that the flow is nonsingular
(no fixed points), the simplest such domain is a solid torus 
$D^2\times S^1$. Invariant solid tori are ubiquitous in volume-preserving
flows on $\real^3$.

In this paper, we extend the programme of Part I to the case of a
solid torus:
\begin{proof}[{\bf Theorem~\ref{mainfluid}}]
{\em	Any steady nonsingular $C^\omega$ Euler flow on $S^1\times D^2$ 
	leaving the boundary invariant possesses a closed flow line.}
\renewcommand{\qed}{}
\end{proof}
\renewcommand{\qed}{\qedsymbol}
Note that this result is independent of the geometry of the solid 
torus. 
For the proof, 
we rely on a dichotomy for $C^\omega$ steady nonsingular Euler flows 
due to Arnold which presents two scenarios: the {\em integrable} and 
the {\em Beltrami} fields. In part I of this series, it was shown that 
the crucial subclass of Beltrami fields is in fact equivalent to the 
class of Reeb fields associated to contact forms (see \S\ref{belreeb}). 
Thus, after dealing with the integrable cases, we reduce the problem of 
periodic orbits for Euler fields to the problem of periodic orbits 
for Reeb fields, \ie, the Weinstein conjecture:

\begin{proof}[{\bf Theorem~\ref{maincontact}}]
{\em 	Every Reeb field on $S^1\times D^2$ tangent to the
	boundary possesses a periodic orbit.}
\renewcommand{\qed}{}
\end{proof}
\renewcommand{\qed}{\qedsymbol}

The work of Hofer \cite{Hof93} resolved the Weinstein conjecture
on the 3-sphere, and on manifolds with nontrivial $\pi_2$ by considering
sequences of pseudo-holomorphic curves in symplectizations. The 
heart of our proof for the solid torus likewise lies in Hofer's results
on pseudo-holomorphic curves. A simple generalization is impossible, 
however, for the reason that a certain class of Reeb fields (namely
those associated to {\em tight} structures --- see \S\ref{contact} --- 
requires different techniques on different manifolds. The key 
steps in the proof of \tm{maincontact} are a careful analysis of
virtually overtwisted contact structures on the solid torus, along 
with an application of the recent work of Hofer et al. on finite
energy foliations.

In \S\ref{fluid} we recall the relevant definitions from hydrodynamics and
discuss the relations with contact geometry. This immediately allows us
to reduce \tm{mainfluid} to the Weinstein conjecture.  
In \S\ref{contact} we turn to contact geometry, 
collecting the standard facts we need. A characterization of 
virtually overtwisted contact structures on the solid torus is obtained in 
\S\ref{virtuallyot}.  In \S\ref{hofer} we discuss Hofer's approach to
the Weinstein conjecture. In particular, we recall the use of 
pseudo-holomorphic curves in contact geometry. The Weinstein conjecture 
on solid tori is proved in \S\ref{Weinstein}.

\section{Fluid flows on solid tori}\label{fluid}

The dynamical properties of incompressible, inviscid fluid flows are 
described by the Euler equations. For an overview of geometric hydrodynamics 
on Riemannian manifolds, the reader is referred to \cite{AK97,AK92}. 
To describe the Euler equations on a 3-manifold $M$ we must first 
fix a Riemannian metric $g$ and a volume form $\mu$. Note that, 
following \cite{AK97}, we do not assume that $\mu$ is the volume form 
associated to $g$, though we do of course allow that possibility. 
The {\em Euler equations} governing the velocity field $u(t)$ of a 
perfect incompressible fluid may be written as follows:
\begin{equation}\label{eq_Euler}
\begin{array}{rcl}
\frac{\displaystyle (\del\iota_u g)}{\displaystyle \del t} + 
	\iota_u\iota_{\nabla\times u}\mu&=& -dP \\ 
\mathcal{L}_u\mu &=& 0 .
\end{array}
\end{equation}
Here $P(t):M\to\R$ is some time-dependent (Bernoulli) function, 
$\mathcal{L}$ is the Lie derivative, $\iota$ denotes contraction, and 
$\nabla\times u$ is the vorticity, defined by the relation 
$\iota_{\nabla\times u}\mu=d\iota_u g$. We call $u$ an Euler field 
and its flow an Euler flow if $u$ satisfies Equation~\ref{eq_Euler}
for some function $P$. In this paper we concern ourselves only 
with steady nonsingular solutions $u$ to Equation~\ref{eq_Euler}; that is, 
solutions for which all time derivatives vanish and all velocities 
are nonzero. For such Euler fields there exists a powerful dichotomy.

\begin{thm}[Arnold \cite{Arn66}]\label{thm_Arnold}
Let $u$ be a $C^\omega$ nonsingular steady Euler field on a $C^\omega$ 
Riemannian three-manifold $M$. If $\partial M\not=\emptyset$ then assume 
$u$ is tangent to the boundary of $M$. If $u$ is not everywhere colinear 
with its curl, then it has a stratified integral: \ie, there exists a compact 
stratified subset $\Sigma\subset M$ of codimension at least one which 
splits $M$ into a finite collection of cells diffeomorphic to $T^2\times\real$. 
Each $T^2\times\{c\}$ is an invariant set for $u$ having flow conjugate 
to linear flow.  
\end{thm}
This theorem motivates the following classical definitions.
\begin{dfn}\label{def_Beltrami} 
A volume-preserving vector field $X$ on a Riemannian manifold 
$M^3$ is a {\em Beltrami field} if 
it is parallel to its curl: \ie, $\nabla\times X = f X$ for some function 
$f:M\to\real$. We say that a Beltrami field is {\em rotational} if $f$ is 
nowhere zero and that it is {\em irrotational} if $f$ vanishes identically.
\end{dfn}
Furthermore, since the function $f$ above is an integral for $X$,
one may reduce the study of real-analytic nonsingular steady Euler fields
to the following three cases:
\begin{enumerate}
\item
	rotational Beltrami fields;
\item
	irrotational Beltrami fields; and
\item
	stratified integrable fields.
\end{enumerate}
This roughly outlines the plan of the proof of the main theorems.

\subsection{Reeb fields associated to Beltrami flows}\label{belreeb}

Given a rotational Beltrami vector field one may consider the 1-form 
$\alpha:=\iota_X g$.
Since $d\alpha=d\iota_X g=f\iota_X\mu$ and the kernel of $\alpha$ is 
orthogonal to $X$ it is clear that 
\begin{equation}\label{eq_reebfield}
	\alpha\wedge d\alpha\not=0.
\end{equation}
Thus, by definition (see \S\ref{contact}), the plane field 
$\xi=\hbox{ker }\alpha$ is a \df{contact structure}. Associated to $\alpha$ 
is a special vector field, the \df{Reeb field} $X_\alpha$,
defined by the conditions
\begin{equation}
	\iota_{X_\alpha}\alpha=1, \mathspace  \iota_{X_\alpha}d\alpha=0.
\end{equation}
As $\iota_X d\alpha=\iota_X d\iota_X g=\iota_X\iota_X\mu=0$, one 
concludes that $X=hX_\alpha$ for some nonzero function $h:M\to\real$.  
Thus any rotational Beltrami field is a nonzero rescaling of some Reeb 
field (sometimes referred to as a \df{Reeb-like} field). This is the 
easy half of the following:
\begin{thm}[Etnyre \& Ghrist \cite{EG:I}]
\label{thm_Correspondence}
	On a fixed 3-manifold $M$, the class of vector fields which are 
	nonsingular rotational Beltrami fields for some Riemannian 
	structure is equivalent to the class of vector fields which 
	are nonsingular rescalings of the Reeb field of some contact form.
\end{thm}

If $X$ is an irrotational Beltrami field, then the 1-form $\alpha
= \iota_Xg$ still defines a plane field 
$\xi:=\ker\alpha$ transversal to $X$; however, in this case 
$d\alpha=d\iota_X g=\iota_{\nabla\times X}\mu\equiv 0$. The Frobenius 
Theorem implies that $\xi$ generates a codimension-one foliation 
of $M$ transverse to the boundary.

\subsection{Periodic orbits in Euler flows}

In this section we set the stage for the proof of the main theorem:

\begin{thm}\label{mainfluid}
	Any steady nonsingular $C^\omega$ Euler field on $S^1\times D^2$ 
	tangent to the boundary possesses a closed flowline.
\end{thm}
\pf 
Since $u$ is a real-analytic steady Euler field,
Theorem~\ref{thm_Arnold} and the discussion above imply that one 
needs to consider three cases: (1) $u$ is a Reeb vector field for 
a contact structure; (2) $u$ preserves a transverse foliation; or (3) 
$u$ has a stratified integral. Case (1) is precisely the Weinstein 
Conjecture on the solid torus, which we prove
in \S\ref{Weinstein}. A theorem of Tischler \cite{Tis70} 
implies that, in case (2), the foliation is actually by fibers of 
a fibration of $S^1\times D^2$ over $S^1$. Using the exact sequence 
for homotopy groups of a fibration one easily concludes that 
the fiber must be $D^2$. The vector field $u$ is transverse to 
the fibers; thus, any fiber will provide a section to the flow, and 
one concludes the theorem via the Brouwer fixed point 
theorem.
 
We are left to consider an Euler field $u$ with a stratified integral. The
argument for this case essentially mirrors the argument in \cite{EG:I} 
for the integrable case on a Riemannian $S^3$. 
Specifically, the real-analytic
codimension-one (or greater) set $\Sigma$ from Theorem~\ref{thm_Arnold} 
possesses a certain [Whitney] stratification, each stratum of which 
is invariant under the flow of $u$. Thus $\Sigma$ has no zero-strata. 
The collection of essential one-dimensional strata must be nonempty, 
otherwise we would have foliated $D^2\times S^1$ by copies of $T^2$: 
impossible. These one-strata thus 
consist of closed 1-manifolds invariant under the flow.
\qed

\begin{remark}
	In the contact case of the above theorem, one ``frequently'' finds 
	contractible periodic orbits which do not exist in the foliation
	case and do not necessarily exist in the integrable case. Following
	the outline of the above proof combined with the proof of the Weinstein
	conjecture on solid tori below one may extract a computable invariant
	of Euler fields which 
	can be used to detect contractible periodic orbits.
\end{remark}

\section{Contact structures and characteristic foliations}\label{contact}

Recall a {\em contact structure} $\xi$ on a 3-manifold $M$ is 
a plane field in $TM$ that is maximally nonintegrable. For the 
remainder of this work, we may assume $\xi$ to be transversally 
orientable so that $\xi=\hbox{ker }\alpha$ for a nondegenerate 
1-form $\alpha$ (since all the contact structures we encounter 
come with a Reeb field by Theorem~\ref{thm_Correspondence}). Such an 
$\alpha$ is a {\em contact form}.
One can express the nonintegrability of $\xi$ by 
the condition $\alpha\wedge d\alpha\not=0$.
Two contact structures are {\em contactomorphic} if there is a diffeomorphism
of $M$ that takes one of the contact structures to the other.
For proofs of some of the standard results listed below see \cite{Aeb94}.

Given a surface $\Sigma$ in $M$, the contact structure $\xi$ induces a 
singular foliation $\Sigma_\xi$ on $\Sigma$, generated by the 
line field $T\Sigma\cap \xi$, with singularities occurring at points 
where $T_p\Sigma=\xi_p$. This is known as the {\em characteristic 
foliation}. 

\begin{lem}[Moser-Weinstein]\label{foldet}
	Two contact structures that induce the same characteristic foliation
	on a surface are contactomorphic in a neighborhood of the surface.
\end{lem}

A contact structure is \df{overtwisted} if there exists an
embedded disc in $M$ whose characteristic foliation contains a limit 
cycle, otherwise it is called \df{tight}. If $\xi$ is tight and there is a 
finite cover of $(M,\xi)$ that is overtwisted then $\xi$ is
called \df{virtually overtwisted}. 

Generically,  
the singularities in $\Sigma_\xi$ are either of elliptic or hyperbolic type.
Moreover, each singularity of $\Sigma_\xi$ is assigned a sign depending
on whether or not the orientations on the plane field $\xi$ and $T\Sigma$ 
agree or not. Of paramount importance in detecting overtwisted 
discs is the Elimination Lemma:

\begin{lem}[Elimination Lemma \cite{Gir91,Eli92}]\label{elim}
	Let $\Sigma$ be a surface in a contact 3-manifold $(M,\xi)$.
	Assume that $p$ and $q$ are singular points in $\Sigma_\xi$ 
	which are of different type (one elliptic, one hyperbolic)
	yet have the same sign. Finally assume that there exists
	a leaf $\gamma$ in $\Sigma_\xi$ that connects $p$ to $q$.
	Then, given any small neighborhood $U$ of $\gamma$ in $M$, there 
	exists an arbitrarily $C^0$-small isotopy of $\Sigma$, fixed outside
	of $U$, which removes all singularities of $\Sigma_\xi$ within $U$.
\end{lem}

This lemma is most effective when used in conjunction with:

\begin{lem}[\cite{Etn96,EF98}]\label{monodromy}
	Let $\gamma$ and $\gamma'$ be two leaves in the characteristic
	foliation $\Sigma_\xi$ both ending in an elliptic singularity
	$e$.  By a $C^0$-small isotopy of $\Sigma$ near $e$ we may arrange
	for $\gamma\cup\gamma'\cup\{e\}$ to be a smooth curve in the new 
	characteristic foliation. 
\end{lem}

The above two lemmas allow us to cancel singularities of the same sign so 
that the leaf joining the elliptic point to the hyperbolic point and
any other leaf touching the elliptic point join to form
a smooth leaf after the cancelation. For example, the two branches
of the unstable (or stable, depending on the orientation) manifold
of a hyperbolic point terminating in an elliptic fixed point of the same 
sign will form a smooth leaf after the cancelation
of the positive (negative respectively) singularities.

The next result tells us how to smooth corners of a surface without 
changing the characteristic foliation.

\begin{lem}[Makar-Limanov \cite{Mak98}]\label{smooth}
	Suppose $\Sigma$ and $\Sigma'$ are two surfaces with boundary
	in $(M^3,\xi)$ that intersect transversally along their boundaries.
	If $\gamma:=\Sigma\cap \Sigma'$ is transverse to $\xi$ then we may
	isotope $\Sigma$ and $\Sigma'$ in an arbitrary neighborhood of 
	$\gamma$ so that $\Sigma\cup\Sigma'$ is a smooth surface and 
	the isomorphism type of the characteristic foliation is 
	unchanged throughout the isotopy.
\end{lem}

If $\gamma$ is a curve embedded in a contact manifold $(M,\xi)$ 
we say $\gamma$ is a \df{transversal curve} if $\gamma$ is transversal 
to $\xi$ at each point of $\gamma$.  If $\gamma$ is nullhomologous
then there is an embedded surface $\Sigma$ such that 
$\partial \Sigma=\gamma$.  Note that $\xi\vert_{\Sigma}$ is a trivial bundle 
so we may choose a nonzero section $s$ of $\xi\vert_{\Sigma}$. 
Using this section to push off a copy $\gamma'$ of $\gamma$ leads 
naturally to the definition of the \df{self-linking number of $\gamma$}
as the intersection number of $\gamma'$ with $\Sigma$: 
\begin{equation}
	\ell(\gamma;\Sigma) := \gamma'\cdot \Sigma
\end{equation}
In subsequent sections, we will consider transversal meridians of 
the boundary of contact solid tori. Proofs of the main theorems 
differ depending on the self-linking numbers of the meridians.
For some of the many properties and applications of self-linking 
numbers in contact geometry, see \cite{Aeb94,Ben83,ET97}.

\section{Covers of tight contact tubes}\label{virtuallyot}

Recall the definitions of tight and overtwisted from \S\ref{contact}:
overtwisted contact structures possess discs with a limit cycle in 
the characteristic foliation. 
Let $\xi$ denote a tight contact structure on a solid torus 
$V=S^1\times D^2$.  We further 
assume that the characteristic foliation $(\del V)_\xi$ 
on the boundary of $V$ is nonsingular and contains no
Reeb components (\ie, annuli possessing no transversals from one 
boundary component to the other). 
In this situation we may choose a meridional curve 
$\mu$ on $\del V$ that is positively transverse
to $(\del V)_\xi$. Denote by $m_\xi$ the self-linking number 
of $\mu$ relative to the
disc $D$ that $\mu$ bounds.  Given that $\xi$ is tight, one knows 
from the inequalities in, \eg, \cite{Eli92}, that:
\begin{equation}
m_\xi\leq -\chi(D)=-1 .
\end{equation}
We may now state:

\begin{thm}\label{otcover}
	Let $\xi$ be a tight contact structure on the solid torus $V$
	for which $(\partial V)_\xi$ is a linear foliation (conjugate
	to  a foliation by lines of constant slope).
	The contact structure $\xi$ is virtually overtwisted if and only if 
	the self-linking number of the meridian is less than $-1$.
\end{thm}

This theorem is also true in the case where the foliation is not
linear; however, we provide the [much shorter] proof of the simpler
case as this is all we require for the sequel.
We begin with some preliminary steps.
In \cite{Eli92} it was shown that all of the negative 
elliptic and positive hyperbolic
singularities of a characteristic foliation on a disc 
may be eliminated by a $C^0$-small isotopy of the disc fixed near the
boundary. After this elimination procedure on our meridional disc,  
we are left with $e_+$ positive elliptic and $h_-$ negative hyperbolic 
singularities.  One may easily check that
\begin{equation}
\label{eq_ellhyp}
m_\xi=-(e_+ + h_-).
\end{equation}

Note that since $(\del V)_\xi$ is nonsingular we may use 
this foliation to define a return map $\Phi$ on $\mu$. Recall, 
associated to $\Phi$ (or to any circle diffeomorphism) is 
a rotation number $r(\Phi)\in[0,1)$. 

\begin{prop}\label{prelim}
	Suppose $V$ deformation retracts to a solid torus $V'\subset V$ 
	for which $(\del V')_\xi$ is nonsingular. If, for any such $V'$, 
	the return map $\Psi'$ for the meridian $\mu'$ of $V'$ has an 
	irrational rotation number and the self-linking number of $\mu'$ 
	is strictly less than $-1$, then $V$ has an overtwisted finite 
	cover.
\end{prop}

\begin{figure}[ht]
	{\epsfxsize=5.0in\centerline{\epsfbox{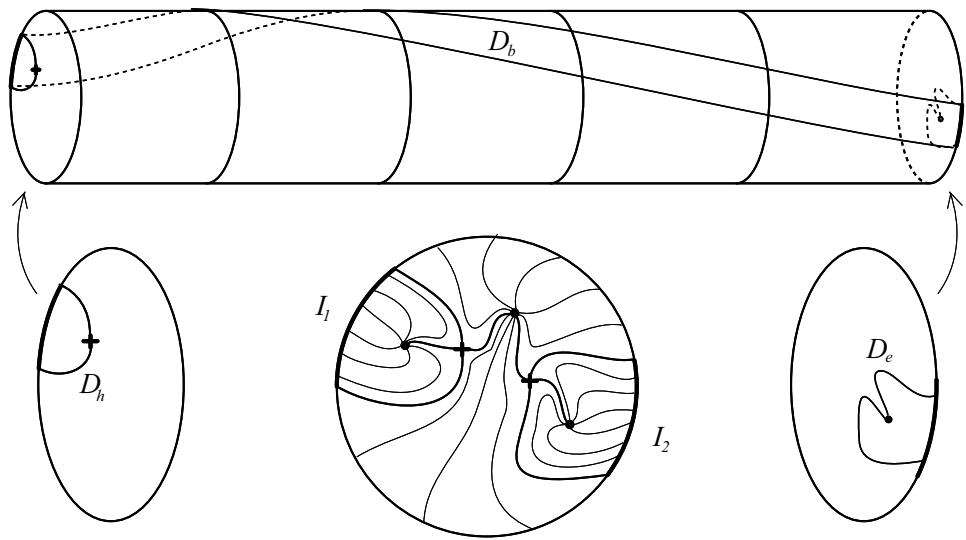 }}}
	\caption{Above, the $(n+1)$-fold cover of $V'$ ($n$ fundamental domains
	shown, above) contains an overtwisted disc. Below, center, the 
	characteristic foliation $D'_\xi$. }
	\label{fig:otcover}
\end{figure}

\pf
Since the self-linking of $\mu'$ is strictly less than 
$-1$ there will be at least one hyperbolic point, from 
Equation~\ref{eq_ellhyp}. If there is only one such point, 
$\mu'$ is divided into two open intervals $I_1$ and $I_2$ by 
the ends of the unstable manifold of the hyperbolic point.  
Otherwise we can find a pair of hyperbolic points, $h_1$ and $h_2$, 
each of whose unstable manifold $W^u(h_j)$ divides the meridional 
disc $D'$ into two subdiscs, one of which encloses a unique 
singular point of elliptic type. Denote by $I_j\subset\mu'$ 
the arc subtended by $W^u(h_j)$ which, together with $W^u(h_j)$, 
bounds the subdisc containing a unique elliptic point (as in 
Figure~\ref{fig:otcover}). 

We claim that some iterate of $\Psi'$ maps $I_1$ into $I_2$ (or 
vice-versa). Indeed since $\Psi'$ is topologically conjugate to an 
irrational rigid rotation, we may argue as if it is such a rotation. 
Iterates of $\Psi'$ map the clockwise endpoint of $I_1$
arbitrarily close to the clockwise endpoint of $I_2$. So either
an iterate of $\Psi'$ maps $I_1$ into $I_2$, or $I_2$ 
into $I_1$, or the $I_j$'s have exactly the same length. 
If they all have the same length, then a small perturbation of $D'$ 
near one of the hyperbolic points will change
the length of one of the $I_j$, thus proving our claim.

Suppose $(\Psi')^n$ maps $I_1$ into $I_2$. Then the $(n+1)$-fold 
cover $\tilde{V}$ of $V$ contains an overtwisted disc as illustrated 
in Figure~\ref{fig:otcover}. The cover is composed of $n+1$ copies 
of $V$ cut along its meridional disc, inside of which there are
$n+1$ copies of $V'$, labeled $V'_i$, cut along $D'$.  In $V'_1$, say,
one has a subdisc $D_h$ of $D'$ cut out by $W^u(h_2)$ and $I_1$. 
Let $D_b$ be the disc consisting of all leaves in the characteristic
foliation  of $\del\tilde{V'}$ emanating from $I_1$ in $\del V'_1$ and 
terminating within $I_2$ in $\del V'_n$.  Finally  let $D_e$ be the
subdisc of $D'$ in $V_n$ consisting of leaves of the characteristic foliation
of $D'$ emanating form the interval $D_b\cap I_2$ union the elliptic point
to which they limit. We may now use Lemma~\ref{smooth} to smooth the corners
of $D_h\cup D_b\cup D_e$ and obtain a disc 
$D_o$ without changing the characteristic foliation. 
So $\partial D_o$ is tangent to the characteristic foliation and contains
exactly one elliptic and one hyperbolic singularity of the same sign. 
Using the Elimination Lemma~\ref{elim} and Lemma~\ref{monodromy} we may 
cancel these singularities leaving $\partial D_o$ a closed
leaf in the characteristic foliation: an overtwisted disc.
\qed

\begin{prop}\label{irratorus}
	There is a near-identity deformation retraction of $V$ 
	to $V'$ such that $(\del V')_\xi$ is nonsingular and irrational. 
\end{prop}
\pf 
Recall that we work under the assumption that $(\del V)_\xi$ is 
a linear foliation. As such, it is a simple matter to construct a 
solid torus $W$ in $(\R^2\times S^1,\xi_0)$, where $\xi_0=\ker(d\phi+r^2d\theta)$ 
[in polar coordinates] so that $(\partial W)_{\xi_0}$ and 
$(\partial V)_\xi$ agree: simply choose the tube $\{r\leq\sqrt{-\kappa}\}$,
where $\kappa<0$ is the slope of $(\del V)_\xi$. 
Using Lemma~\ref{foldet}, a neighborhood of $\del V$ in $V$ is 
contactomorphic to a neighborhood of $\del W$ in $W$ ---
note that the contactomorphism cannot be extended over all of $W$ 
since the meridians of $\del V$ and $\del W$ have different self-linking numbers. 
Shrinking $\del W$ radially yields a 1-parameter family of nearby 
tori with linear foliations varying nontrivially and continuously. 
Thus, tori with irrational foliations exist near $\del W$. Pulling back 
the deformation retraction of $W$ by the contactomorphism yields 
the desired solid torus $V'\subset V$. 
\qed

\pf {\em (of Theorem~\ref{otcover})}
Given $V$ with $m<-1$ we may deformation retract this slightly to obtain
the solid torus $V'$ such that $(\del V')_\xi$ is linear and irrational.
The restriction of $D$ to $V'$ is still a meridional
disc for $V'$ with transverse boundary having self-linking number $m$. 
Proposition~\ref{prelim} then implies an overtwisted cover.

In the case where $m=-1$, a result of Makar-Limanov 
\cite{Mak98} states that $(V,\xi)$ is contactomorphic
to $W=\{(r,\theta, \phi) \in \R^2\times S^1\vert r\leq 
f(\theta, \phi)\}$ for some positive function $f:T^2\to \R$ with the 
tight contact structure, $d\phi+r^2 d\theta$, on $W$.  
By lifting along the $\phi$ coordinate, one obtains 
$\tilde{W}$ an infinite cylinder in $(\R^2\times\R^1 , 
d\tilde{\phi} + r^2d\theta)$. This is the standard
tight contact structure on $\real^3$ (in polar coordiantes).
Since pulling back the contact form on $W$ 
(and hence on $V$) to its universal cover yields a tight
structure, the same result holds for all finite covers.
\qed

\section{Pseduo-holomorphic curves and periodic orbits}\label{hofer}

\subsection{Reeb fields and the Weinstein conjecture}

Recall the definition of a {\em Reeb field}, $X_\alpha$, associated to a 
contact form $\alpha$:
\begin{equation}
	\iota_{X_\alpha}\alpha=1, \mathspace  \iota_{X_\alpha}d\alpha=0.
\end{equation}
Certain questions in Hamiltonian dynamics can be reformulated in 
terms of the dynamics of $X_\alpha$.  This relation and results of 
Rabinowitz \cite{Rab79} and Weinstein \cite{Wei79} concerning periodic 
orbits in Hamiltonian dynamics led Weinstein to pose the following:

\begin{conj}[Weinstein Conjecture]
	For each contact form on a closed 3-manifold 
	the corresponding Reeb vector field has
	a periodic orbit.
\end{conj}

Hofer \cite{Hof93} has recently made extraordinary
progress on the Weinstein conjecture.  Among other things he has shown
the following:

\begin{thm}\label{perorbit}
	Any contact form $\alpha$ associated to an overtwisted 
	contact structure on a closed 3-manifold $M$ possesses a closed
	orbit in its Reeb field which is of finite (perhaps trivial) 
	order in $\pi_1(M)$.
\end{thm}

It has been observed by some experts that Hofer's theorem is still 
true for manifolds with boundary provided the Reeb vector field is 
tangent to the boundary. Unfortunately neither this result nor its 
proof have appeared in the literature.  So for the sake of completeness
we briefly sketch the proof of \tm{perorbit} noting the necessary 
modifications to make it valid for manifolds with boundary.

\subsection{A review of pseduo-holomorphic curves}\label{J}

The main tool in the proof of \tm{perorbit} 
is the use of pseudo-holomorphic 
curves in the symplectization of the contact manifold.  
Given a contact 3-manifold $(M,\xi)$ with defining form $\alpha$ 
there is an induced symplectic form 
	\begin{equation}
	\omega=d(e^t\alpha)=e^t(dt\wedge\alpha+d\alpha)
	\end{equation}
on $W=\R\times M$. Choose a complex structure $J_\xi:\xi\to \xi$ on $\xi$ 
so that $\alpha(v,J_\xi v)>0$ for all $v\in \xi$, then define an almost
complex structure $J$ on $W$ by 
	\begin{equation}\label{complex-structure}
	J(a,b)(h,k)=(-\alpha_b(k),J_\xi(b)\pi(k)+hX_{\alpha}(b)),
	\end{equation}
where $(h,k)\in T_{(a,b)}(\R\times M)$, $X_\alpha$ is the Reeb vector 
field for $\alpha$, and $\pi:TM\to\xi$ is projection to $\xi$ along 
$X_\alpha$. Now if $(S,j)$ is a closed Riemannian surface
and $\Gamma$ is a finite subset of $S$ then a map $u:S-\Gamma\to W$ 
is called \df{$J$-holomorphic} (or, \df{pseudo-holomorphic}, if no $J$ is 
specified) if
	\begin{equation}
	du\circ j=J\circ du .
	\end{equation}
One may readily check that that if $\Gamma=\emptyset$ then $u$ is constant.
Given a map  $u:S-\Gamma\to W$ one defines the \df{energy} to be
	\begin{equation}
	E(u)=\sup_{\phi\in\Sigma} \int_{S-\Gamma} 
		u^* d(\phi \alpha),
	\end{equation}
where $\Sigma$ is the set of all smooth maps $\phi:\R\to [0,1]$ 
satisfying $\phi'\geq 0$ and $\phi\alpha$ is the 1-form defined by 
$(\phi\alpha)(a,b)(h,k)=\phi(a)\alpha_b(k)$. In \cite{Hof93} it was shown
\begin{thm}\label{finite-energy-surface}
	If there is a finite energy 
	nonconstant $J$-holomorphic map $u:S- \Gamma\to W$, 
	then $X_\alpha$ has a periodic orbit.
\end{thm}
This theorem is proved by examining the behavior of $u$ near the punctures 
in $S-\Gamma$.

Thus to prove \tm{perorbit} we need merely find a finite energy 
nonconstant $J$-holomorphic map into $W$ when $\xi$ is overtwisted. 
To this end consider an overtwisted disc $\mathcal{D}$ in $M$.
Orient $\mathcal{D}$ so that the unique elliptic point $e$ in 
the characteristic foliation is positive (\ie,  $\xi_e$ defines 
the orientation on $\mathcal{D}$). One then uses the following:
\begin{thm}[Bishop]
	There is a continuous map 
	$$\Psi:D\times [0,\epsilon)\to W$$
	so that for each $u_t=\Psi(\cdot,t)$
	\begin{itemize}
		\item $u_t:D\to W$ is $J$-holomorphic.
		\item $u_t(\partial D)\subset (\mathcal{D}- \{e\})
			\subset \{0\}\times M$.
		\item $u_t\vert_{\partial D}:\partial D\to 
			(\mathcal{D}-\{e\})$ has winding number 1.
		\item $\Psi\vert_{D\times(0,\epsilon)}$ is a smooth map.
		\item $\Psi(z,0)=e$ for all $z\in D$.
	\end{itemize}
\end{thm}
The map $\Psi$ is called a \df{Bishop filling}.
Using an implicit function theorem Hofer finds a \df{maximal Bishop filling}
$\Psi_{\hbox{max}}:D\times[0,1)\to W$. It is important to note that 
$\Psi_{\hbox{max}}(\del D\times[0,1))$
cannot fill all of $\mathcal{D}$, which can be deduced from 
the result \cite{Hof93} that the map $u_t\vert_{\partial D}:\partial D\to 
\mathcal{D}$ is an embedding which is transversal to the characteristic 
foliation on $\mathcal{D}$.
One may then argue that there is a sequence of 
$t_k\to 1$ and $z_k\to z_0$ so that 
$|\nabla\Psi_{\hbox{max}}(z_k,t_k)|\to\infty$,
since, if this were not the case, the sequence $\Phi_{\hbox{max}}(\cdot,t_k)$ 
would converge to a $J$-holomorphic disc which would allow us to extend 
$\Psi_{\hbox{max}}$ to a larger Bishop filling. After reparameterization,
one may assume that the $z_k$ are bounded away from the boundary of $D$: 
thus the gradients are blowing up on the interior of $D$.

Following \cite{Hof93},
assume that $z_k=0$ for all $k$ and the norm of the gradient of 
$\Psi_{\hbox{max}}(\cdot,t_k)$ is maximal at $0$. Define a sequence of 
maps $v_k:D_k\to W$ where $D_k$ is a disc of radius 
$R_k=|\nabla\Psi_{\hbox{max}}(0,t_k)|$ by
	\begin{equation}
	v_k(z)=(a_k(z/R_k)-a_k(0), u_k(z/R_k)),
	\end{equation}
where $\Psi_{\hbox{max}}(z,t_k)=(a_k(z),u_k(z))$. One may then check 
that the $v_k$ converge to a nonconstant $J$-holomorphic finite 
energy plane $v:\C\to W$. Thus \tm{finite-energy-surface} implies 
$X_\alpha$ has a periodic orbit (since $\C=S^2-\{\infty\}$).

The proof of \tm{perorbit} is very implicit in the sense that one 
has no information on the placement of the periodic orbit given 
the placement of the overtwisted disc. In order to extend the proof
to manifolds with boundary, we will consider how the Bishop filling 
interacts with the boundary of the manifold: we require the following 
definition.
\begin{dfn}
Let $(W,J)$ be a 4-dimensional almost complex manifold.  If $M$ is a 
3-dimensional submanifold then there exists a unique hyperplane field 
of complex tangencies in $TM$.  By this we mean that there is a 
2-dimensional subbundle $\eta$ of $TM$ such that $\eta$ is 
$J$-invariant (and hence $J$ is a complex structure on 
$\eta$). Choose a defining 1-form
$\beta$ on $M$ such that $\eta=\ker\beta$. The \df{Levi form}, $L$, 
is defined to be the restriction of $d\beta(\cdot,J\cdot)$ to $\eta$.  
If $L$ is identically zero we say that $M$ is \df{Levi flat} 
(this implies that $\eta$ defines a codimension-one foliation of $M$). 
\end{dfn}
One may now easily verify:
\begin{lem}
	Let $\alpha$ be a contact 1-form on $M^3$ whose Reeb vector 
	field $X_\alpha$ is tangent to $\partial M$. Then the boundary 
	of $W=\R\times M$ is Levi-flat with respect to the almost complex 
	structure $J$ in \eqn{complex-structure}; more specifically, 
	$\partial W$ is foliated by the complex surfaces $\R\times \gamma$ 
	where $\gamma\subset\partial M$ is an orbit of $X_\alpha$.
\end{lem}

Finally we need to recall how $J$-holomorphic curves intersect.
\begin{thm}[McDuff \cite{Mcd91}]\label{thm_McDuff}
	Two closed distinct $J$-holomorphic curves $C$ and $C'$ in an 
	almost complex 4-manifold $(W,J)$ have only a finite number 
	of intersection points. Each such point contributes
	a positive number to the algebraic intersection number $C\cdot C'$. 
\end{thm}

We are now ready to prove:
\begin{thm}\label{perorbit-with-boundary}
	If $\alpha$ is associated to an overtwisted 
	contact structure on a compact 3-manifold with boundary
	and the Reeb vector field is tangent to the boundary, then the
	Reeb vector field has a	closed orbit which is of 
	finite (perhaps trivial) order in $\pi_1(M)$.
\end{thm}
\pf
We begin by completing $M$ to a closed manifold $M'$ containing $M$ and 
extending $\alpha$ to $\alpha'$, a contact 1-form over $M'$, in
the standard manner. Let $W'$ be the associated 
symplectization of $M'$. From the outline of \tm{perorbit} above, 
there exists a maximal Bishop family of
$J$-holomorphic discs $\Psi:D^2\times[0,1)\to W'$ for some standard 
overtwisted disc $\mathcal{D}$ in $M$. Note that 
$\Psi(\partial D,t)\subset \mathcal{D}\subset \{0\}\times M$. We now claim 
that $\Psi(D,t)\subset\real\times M$ as well.  
Indeed, if this is not the case, then one of the $D_t=\Psi(D,t)$ 
would touch $\partial W$ tangentially and thus, since $\partial W$ is 
Levi-flat, $D_t$ would intersect the $J$-holomorphic curve 
$C_\gamma=\R\times\gamma$ (where $\gamma$ is the orbit of the 
Reeb flow on $\del M$ passing through the point of intersection). This 
contradicts \tm{thm_McDuff} since the algebraic intersection of $D_t$ 
and $C_\gamma$ is zero.

Recall one obtains a finite energy plane $v:\C\to W'$ by rescaling the 
Bishop family near the points where the gradient is blowing up.
But since all the $\Psi(D,t)$ lie in $M$, so does $v(\C)$, implying that one
has a periodic orbit of $X_{\alpha'}$ within $M$: a periodic 
orbit of $X_\alpha$.
\qed

\subsection{Finite energy foliations and surfaces of sections}

In this section we discuss some recent work of Hofer et al. 
\cite{Hof98,HWZ99}  
that will be needed to complete the proof of \tm{maincontact}. 
\begin{dfn}
	Let $\alpha$ be a contact form on a 3-manifold $M$ and $J_\alpha$ 
	be a complex structure on $\xi=\hbox{ker}(\alpha)$ as in the 
	beginning of the previous section. A \df{finite energy foliation} 
	of $M$ is a 2-dimensional foliation $\mathcal{F}$
	of $W=\R\times M$ which is invariant under translation along $\R$ and
	whose leaves are $J$-holomorphic surfaces having uniformly bounded 
	energies.
\end{dfn}
Several useful facts concerning finite energy foliations 
appear in \cite{Hof98,HWZ99}:
\begin{lem}[Hofer et al. \cite{Hof98,HWZ99}]
	Let $\mathcal{F}$ be a finite energy foliation of $(M,\alpha)$, then
	\begin{itemize}
		\item If $F$ is a leaf of $\mathcal{F}$ invariant under 
			some translation then $F=\R\times P$ where $P$ is
			a periodic orbit of $X_\alpha$.
		\item If a leaf $F$ is not invariant under any translation
			then its projection $\hat{F}$ to $M$ is an embedded
			submanifold of $M$ transversal to $X$.
		\item If the projection $\hat{F}$ and $\hat{G}$ of two leaves
			of $\mathcal{F}$ intersect in $M$ then $F$ is a 
			translate of $G$.
	\end{itemize}
\end{lem}
This lemma implies that one obtains a foliation on the complement 
of some periodic orbits in $M$ which is transverse to the Reeb flow.

Finite energy foliations of $S^3$ are in some sense generic.  To make this 
precise, recall that a periodic orbit in a Reeb flow is \df{nondegenerate} 
if the linearized Poincar\'e return map associated to the orbit does 
not have 1 as an eigenvalue. A contact form $\alpha$ on a 3-manifold $M$ 
is called \df{nondegenerate} if all the periodic orbits are nondegenerate.
A result in \cite{HWZ98} asserts that for a fixed contact form $\alpha$ on
$M$, the set of positive functions $f$ such that $f\alpha$ is 
nondegenerate is a Baire set in $C^\infty(M,(0,\infty))$.

Given $\alpha_0$ the standard contact form on $S^3$,
the main theorem concerning the 
existence of finite energy foliations is:
\begin{thm}[\cite{Hof98,HWZ99}]
\label{fef}
	If $\alpha$ is a nondegenerate tight contact form on $S^3$, 
	then there is a Baire set of admissible complex 
	structures $J$ on $\xi$ for which 
	$(S^3,\alpha,J)$ admits a finite energy foliation.
\end{thm}

\section{The weinstein conjecture for solid tori}\label{Weinstein}

In this section we prove the Weinstein conjecture on the solid
torus. Specifically,
\begin{thm}\label{maincontact}
	Every Reeb field on $S^1\times D^2$ tangent to the
	boundary possesses a periodic orbit.
\end{thm}

Of course, this applies to general contact 3-manifolds possessing 
invariant solid tori; however, the hypotheses are, for a general Reeb flow, 
almost as hard to verify as the existence of a closed orbit in the first 
place. Nevertheless, this result is useful in certain specific examples.

\begin{lem}\label{linear}
	If $T$ is a torus invariant under the flow of a Reeb field $X$ 
	(associated to $\xi$) and $X$ has no periodic orbits on $T$, then 
	the characteristic foliation $T_\xi$ is linear (up to conjugacy).
\end{lem}

\pf
Since $X$ is tangent to $T$ the characteristic foliation is nonsingular.
We may assume that the characteristic foliation is not a linear 
foliation by meridional curves. 
In the case where there does not exist a transversal meridian (\ie, 
there is a Reeb component in $T_\xi$), then $X$ must have a periodic 
orbit (since $X$ must point transversely into or out of the Reeb 
component, hence limiting to a periodic orbit).  We may thus pick a 
transversal meridian and consider the rotation number associated 
to $T_\xi$.  If the rotation number is irrational then the characteristic 
foliation is conjugate to a linear foliation. If
the rotation number is rational then there is a closed leaf $L$ in the 
foliation.  Since the flow of $X$ preserves both $\xi$ and $T$, it also 
preserves $T_\xi$. Thus $T_\xi$ is simply a foliation by rational curves.
\qed

\proof[Proof of Theorem~\ref{maincontact}:]
Assume there are no periodic orbits on the boundary.
Suppose $X$ is associated to the contact form $\alpha$ 
and the contact structure $\xi$. If $\xi$ is overtwisted 
then we are done by \tm{perorbit-with-boundary}; thus, assume 
that $\xi$ is tight.  Since $X$ is tangent to the boundary, 
the characteristic foliation on $T$ is nonsingular. We may assume that 
the foliation is linear by Lemma~\ref{linear} , and that it is not 
composed of meridians by tightness. Hence, we may choose  
the meridian $\mu$ transverse to $T_\xi$. 
If the self-linking number $m:=\ell(\mu)$ is 
less than $-1$, \tm{otcover} implies that some finite cover of 
$(S^1\times D^2,\alpha,\xi)$ is overtwisted and thus has a periodic 
orbit in its Reeb flow. Under the covering map, flowlines
are mapped to flowlines and hence $X$ must also have a periodic orbit.

We are left to consider the case when $m=-1$.  Assume that the Reeb 
field does not possess a periodic orbit. Then it follows from 
\cite{Mak98}
that there exists a map $f$ from $S^1\times D^2$ to a neighborhood $V$ 
of a transversal unknot in $S^3$ such that $f_*(\xi)=\xi_0$ where 
$\xi_0$ is the standard tight contact structure on $S^3$. We may thus 
push $\alpha$ forward to $V$ and extend it to a contact form ${\alpha'}$ on 
all of $S^3$. Thus we have a Reeb vector field $X_{\alpha'}$ associated to the 
tight contact structure on $S^3$. All the periodic orbits of $X_{\alpha'}$ 
must by assumption lie outside $V$. 

\begin{lem}
\label{nondegenerate_rel}
	One can perturb ${\alpha'}$ fixing $V$ so that it is 
	a nondegenerate contact form.
\end{lem}
\pf 
This proof is the relative version of \cite[Prop. 6.1]{HWZ98}. Briefly, 
one embeds $(S^3,{\alpha'})$ into the symplectization 
$(\real\times S^3,d(e^t{\alpha'}))$ as $\{0\}\times S^3$, identifying the
Reeb field of ${\alpha'}$ with the induced Hamiltonian field on 
the hypersurface. Perturbing the hypersurface $\{0\}\times S^3$ is 
equivalent to perturbing the Hamiltonian field. The theorem of 
Robinson \cite[Thm. 1.B.iv]{Rob71} states that nondegenerate Hamiltonian 
fields are residual among hypersurfaces. In particular, this result
holds for open manifolds in the strong $C^\infty$ topology on the 
perturbations: thus, perturb the hypersurface on the complement
of $V$, with the perturbation going to zero quickly near the 
boundary. This yields a nondegenerate hypersurface which, since 
contact forms are open in the space of 1-forms, implies a nondegenerate 
contact form on $S^3$ which agrees with ${\alpha'}$ on $V$.  
\qed

Given this, \tm{fef} yields a finite energy foliation of $S^3$.  
Let $\mathcal{F}$ denote the foliation transversal to the Reeb flow on  
the complement of some finite set of periodic orbits in $S^3-V$.  
Since $\partial V$ is invariant under the flow of $X_{\alpha'}$ it is 
transversal to $\mathcal{F}$. Thus $\mathcal{F}\cap \partial V$ is a 
foliation of $\partial V$ by circles. Moreover $\mathcal{F}\cap V$ is 
a foliation of $V$ by either discs or annuli; however, the presence of 
any annuli would clearly contradict the fact that $\mathcal{F}$ intersects 
$\partial V$ transversally. Thus there is a foliation of $V$ by discs 
transversal to the Reeb field. Following the proof in the case
of a stratified integral, this must be a foliation by meridional 
discs. An application of the Brouwer fixed point theorem concludes 
the proof.
\qed

\begin{remark}
We note that an alternate approach to the final step in the proof of 
\tm{maincontact} exists. Instead of using finite-energy foliations, one can
use the (currently developing) {\em contact homology} of Hofer and 
Eliashberg \cite{Eli99}. This is a homology theory that is 
defined in terms of a contact 1-form (actually a
Reeb vector field) but depends only on the underlying contact structure.
The chain groups for this homology are generated by periodic 
orbits in the Reeb flow and the boundary operator is defined using 
pseudo-holomorphic curves in the symplectization that limit to the
periodic orbits in various ways.  It seems the discussion in 
\S\ref{J} is sufficient to allow one to define this contact homology 
for manifolds with boundary (if the implicated Reeb fields are all 
tangent to the boundary).
Then, using standard models for all the 
universally tight contact structures, one can compute that the contact 
homology is non-trivial and thus there must be periodic orbits in any 
Reeb field associated to a universally tight contact form on $S^1\times D^2$.
\end{remark}

\bibliographystyle{plain}


\end{document}